\documentclass[reqno]{amsart}

\usepackage{enumerate}

\theoremstyle{plain}
  \newtheorem{theorem}{Theorem}
  \newtheorem{lemma}{Lemma}
  \newtheorem{proposition}{Proposition}

\newcommand\R{\mathbb R}
\newcommand\Z{\mathbb Z}

\DeclareMathOperator{\pt}{pt}
\DeclareMathOperator{\PSL}{PSL}
\DeclareMathOperator{\ssl}{\mathfrak{sl}}
\DeclareMathOperator{\Diff}{Diff}
\DeclareMathOperator{\id}{id}
\DeclareMathOperator{\img}{img}
\DeclareMathOperator{\supp}{supp}
\DeclareMathOperator{\grad}{grad}
\DeclareMathOperator{\ind}{ind}

\begin{document}

\title{Smooth perfectness for the group of diffeomorphisms}

\author{Stefan Haller}

\address{Stefan Haller,
         Department of Mathematics, University of Vienna,
         Nordbergstra{\ss}e 15, A-1090, Vienna, Austria.}

\email{stefan.haller@univie.ac.at}

\author{Tomasz Rybicki}

\address{Tomasz Rybicki,
         Faculty of applied mathematics, AGH University of Science and Technology,
         Al. Mickiewicza 30, 30-059 Krakow, Poland.}

\email{tomasz@agh.edu.pl}

\author{Josef Teichmann}

\address{Josef Teichmann,
         Department of Mathematics, ETH Z\"urich,
         R\"amistrasse 10, 8092 Z\"urich, Switzerland.}

\email{josef.teichmann@math.ethz.ch}

\thanks{The first author acknowledges the support of the Austrian Science Fund, grant P14195-MAT.
        During the final stage of preparation of this manuscript the first author
        was partially supported by the Austrian Science Fund (FWF), grant
        P19392-N13, and the IK I008-N funded by the University of
        Vienna. The second author was partially supported by the Polish Ministry of Science and Higher
        Education. A  part of this work was written when the
        second author visited ETH in Zurich. He thanks for the invitation
        and hospitality. The third author gratefully acknowledges support from the ETH Foundation.}




\begin{abstract}
  Given a result of Herman, we provide a new elementary proof of the
  fact that the connected component of the group of compactly supported diffeomorphisms is perfect and hence simple.
  Moreover, we show that every diffeomorphism $g$, which is sufficiently close to the identity, can be represented
  as a product of four commutators, $g=[h_1,k_1]\circ\cdots\circ[h_4,k_4]$, where the
  factors $h_i$ and $k_i$ can be chosen to depend smoothly on $g$.
\end{abstract}

\maketitle

\section{Introduction and statement of the result}\label{S:intro}

Let $M$ be a smooth manifold of dimension $n$. By $\Diff^\infty_c(M)$ we will denote
the group of compactly supported diffeomorphisms of $M$. This is a regular Lie group
in the sense of Kriegl--Michor, modelled on the convenient vector space $\mathfrak X_c(M)$
of compactly supported smooth vector fields on $M$, see \cite[Section~43.1]{KM97}. A curve $c\colon \mathbb R\to\Diff^\infty_c(M)$ is smooth iff
the associated map $\hat c\colon \mathbb R\times M\to M$ is smooth and its support satisfies the following condition: for every
compact interval $I\subseteq\mathbb R$ there exists a compact set $ C \subseteq M$ so that $\supp(c(t))\subseteq C $, for all $t\in I$,
see  \cite[Section~42.5]{KM97}.
If $M$ is compact, this smooth structure coincides with the well known Fr\'echet--Lie group structure on $\Diff^\infty(M)$.

Given smooth complete vector fields $X_1,\dotsc,X_N$ on $M$, we consider the map
\begin{align}\label{E:K}
K&\colon \Diff^\infty_c(M)^N\to\Diff^\infty_c(M),\qquad
\\\notag
K(g_1,\dotsc,g_N)&:=[g_1,\exp(X_1)]\circ\cdots\circ[g_N,\exp(X_N)].
\end{align}
Here $\exp(X)$ denotes the flow of a complete vector field $X$ at time $1$, and $[k,h]:=k\circ h\circ k^{-1}\circ h^{-1}$ denotes the commutator of two diffeomorphisms
$k$ and $h$. It is readily checked that $K$ is smooth. Indeed, one only has to observe that $K$ maps smooth curves to smooth curves in view
of the characterization of smooth curves from the previous paragraph, cf.\ \cite[Section~27.2]{KM97}.
Writing $\id$ for the identity in $\Diff^\infty_c(M)$, we clearly have $K(\id,\dotsc,\id)=\id$.

A \emph{smooth local right inverse at the identity for $K$} consists of a $c^\infty$-open neighborhood $\mathcal U$ of
the identity in $\Diff^\infty_c(M)$ together with a smooth map
$$
\sigma=(\sigma_1,\dotsc,\sigma_N)\colon \mathcal U\to\Diff^\infty_c(M)^N
$$
so that $\sigma(\id)=(\id,\dotsc,\id)$ and $K\circ\sigma=\id_{\mathcal U}$. More explicitly, we require that
each $\sigma_i\colon \mathcal U\to\Diff^\infty_c(M)$ is smooth with $\sigma_i(\id)=\id$ and, for all $g\in\mathcal U$,
$$
g=[\sigma_1(g),\exp(X_1)]\circ\cdots\circ[\sigma_N(g),\exp(X_N)].
$$
The $c^\infty$-topology \cite[Section~4]{KM97} is the final topology with respect to all smooth curves, its open subsets
are the natural domains for locally defined smooth maps between infinite dimensional manifolds \cite[Section~27]{KM97}.
For compact $M$ the $c^\infty$-topology on $\Diff^\infty(M)$ coincides with the Whitney $C^\infty$-topology, cf.\ \cite[Theorem~4.11(1)]{KM97}.
In general the $c^\infty$-topology on $\Diff^\infty_c(M)$ is strictly finer than the one induced from the Whitney $C^\infty$-topology,
cf.\ \cite[Section~4.26]{KM97}. The latter coincides with the inductive limit topology $\lim_K\Diff^\infty_K(M)$ where $K$ runs
through all compact subsets of $M$, see \cite[Section~41.13]{KM97}.

The aim of this paper is to establish the following two results.

\begin{theorem}\label{T:main}
Suppose $M$ is a smooth manifold of dimension $n\geq2$.
Then there exist four smooth complete vector fields $X_1,\dotsc,X_4$ on $M$
so that the map $K$, see \eqref{E:K}, admits a smooth local right inverse at the identity, $N=4$.
Moreover, the vector fields $X_i$ may be chosen arbitrarily close to zero
with respect to the strong Whitney $C^0$-topology.
If $M$ admits a proper (circle valued) Morse function whose critical points all have index $0$ or $n$,
then the same statement remains true with three vector fields.
\end{theorem}

Particularly, on the manifolds $M=\R^n,S^n,T^n$, $n\geq2$, or the total space of a compact smooth fiber bundle
$M\to S^1$, three commutators are sufficient. Circle-valued Morse theory was initiated by Novikov \cite{Novikov}, see also Pajitnov's
monograph \cite{Pajitnov} and the references therein.
At the expense of more commutators, it is possible to gain further control on the vector fields.
More precisely, we have:

\begin{theorem}\label{T:mainp}
Suppose $M$ is a smooth manifold of dimension $n\geq2$ and set $N:=6(n+1)$.
Then there exist smooth complete vector fields $X_1,\dotsc,X_N$ on $M$
so that the map $K$, see \eqref{E:K}, admits a smooth local right inverse at the identity.
Moreover, the vector fields $X_i$ may be chosen arbitrarily close to zero
with respect to the strong Whitney $C^\infty$-topology.
\end{theorem}

Either of the two theorems implies that $\Diff^\infty_c(M)_o$, the connected component of the identity,
is a perfect group. This was already proved by Epstein \cite{E84}
using ideas of Mather \cite{M74, M75} who dealt with
the $C^r$-case, $1\leq r<\infty$, $r\neq n+1$. The Epstein--Mather proof is
based on a sophisticated construction, and uses the
Schauder--Tychonov fixed point theorem. The existence of a presentation
$g=[h_1,k_1]\circ\cdots\circ[h_N,k_N]$ is guarantied, but without any further
control on the factors $h_i$ and $k_i$. Rough estimates on the number of necessary factors
are well known, too. In view of a result due to Tsuboi \cite[Theorem~8.1]{T81}, one can assume
$N=4^n(n+1)$, provided $g$ is sufficiently close to the identity. 
More refined estimates for certain classes of manifolds can be found in \cite{T08} and \cite{T09}.
That the factors $h_i$ and $k_i$
can be chosen to depend smoothly on $g$ seems to be folklore as well.

Theorem~\ref{T:main} or \ref{T:mainp} actually implies that the universal covering of $\Diff^\infty_c(M)_o$ 
is a perfect group. This result is known, too, see \cite{T74, M84}.
Thurston's proof is based on a result of Herman for the torus \cite{H71,H73}.

Our proof rests on Herman's result, too,
but is otherwise elementary and different from Thurston's
approach. In fact we only need Herman's result in dimension $1$,
which is a more structured situation, see also \cite{yoccoz95}.

Note that the perfectness of $\Diff^\infty_c(M)_o$ implies that
this group is simple, see \cite{E70}. The methods used in \cite{E70}
are elementary and actually work for a rather large class of
homeomorphism groups.

Let us mention that some analogues of Theorems~\ref{T:main} and
\ref{T:mainp} for the homeomorphism groups in the category of
topological manifolds have been obtained in \cite{R11} by using
completely different arguments.

The remaining part of this note is organized as follows: In
Section~\ref{S:herman} we recall the above mentioned result of
Herman and derive a corollary, see Proposition~\ref{P:herman},
which asserts that the statement of Theorem~\ref{T:mainp} holds
true for the torus, $M=T^n$, with $N=3$. Using the exponential law
we then establish a similar statement for diffeomorphisms on open
subsets of $\mathbb R^n$, see Proposition~\ref{P:exp} in
Section~\ref{S:exp}. This construction allows us to circumvent
Thurston's deformation construction, and at the same time
restricts the approach to dimensions $n\geq2$.
In Section~\ref{S:frag} we formulate and prove a smooth
version of the fragmentation lemma, see Proposition~\ref{P:frag},
and give a proof of Theorem~\ref{T:mainp}.
Finally, in Section~\ref{S:proof} we discuss a technique
to reduce the number of commutators which will eventually
lead to a proof of Theorem~\ref{T:main}.

\section{Herman's theorem revisited}\label{S:herman}

Let $T^n:=\R^n/\Z^n$ denote the torus. For $\lambda\in T^n$ we let $R_\lambda\in\Diff^\infty(T^n)$
denote the corresponding rotation. The main ingredient in the proof of 
Theorems~\ref{T:main} and \ref{T:mainp} is the following result of Herman \cite{H73, H71}.

\begin{theorem}[Herman]\label{T:herman}
There exist $\gamma\in T^n$ so that the smooth map
$$
T^n\times\Diff^\infty(T^n)\to\Diff^\infty(T^n),\qquad(\lambda,g)\mapsto R_\lambda\circ[g,R_\gamma],
$$
admits a smooth local right inverse at the identity.
Moreover, $\gamma$ may be chosen arbitrarily close to the identity in $T^n$.
\end{theorem}

Herman's result is an application of the Nash--Moser inverse
function theorem. When inverting the derivative one is quickly led
to solve the linear equation $Y=X-(R_\gamma)^*X$ for given
$Y\in C^\infty(T^n,\R^n)$. This is accomplished using Fourier
transformation. Here one has to choose $\gamma$ sufficiently irrational
so that tame estimates on the Sobolev norms of $X$ in terms
of the Sobolev norms of $Y$ can be obtained. The corresponding
small denominator problem can be solved due to a number
theoretic result of Khintchine.

Below we will make use of the following corollary of Herman's result:

\begin{proposition}\label{P:herman}
There exist smooth vector fields $X_1,X_2,X_3$ on $T^n$ so that the smooth map
$\Diff^\infty(T^n)^3\to\Diff^\infty(T^n)$,
$$
(g_1,g_2,g_3)\mapsto[g_1,\exp(X_1)]\circ[g_2,\exp(X_2)]\circ[g_3,\exp(X_3)],
$$
admits a smooth local right inverse at the identity. Moreover, the vector fields $X_i$ may be chosen
arbitrarily close to zero with respect to the Whitney $C^\infty$-topology.
\end{proposition}

\begin{proof}
This is an immediate consequence of Theorem~\ref{T:herman} and the following observation, cf.\ \cite{H73}.
The finite dimensional Lie group $\PSL_2(\mathbb R)$ acts effectively on the circle $T^1$ and we have smooth embeddings
$T^1\subseteq\PSL_2(\mathbb R)\subseteq\Diff^\infty(T^1)$. Since $\PSL_2(\mathbb R)$ is a simple Lie group, there exist $Y_1,Y_2\in\ssl_2(\mathbb R)$
so that the smooth map
$$
\PSL_2(\mathbb R)^2\to\PSL_2(\mathbb R),\qquad(g_1,g_2)\mapsto[g_1,\exp(Y_1)][g_2,\exp(Y_2)],
$$
admits a smooth local right inverse at the identity. Moreover, $Y_i$ may be chosen arbitrarily close to $0$ in $\ssl_2(\mathbb R)$.
Taking the product of $n$ copies of such local right inverses and using the smooth embeddings
$T^n\subseteq\PSL_2(\mathbb R)^n\subseteq\Diff^\infty(T^1)^n\subseteq\Diff^\infty(T^n)$,
the statement follows readily from Theorem~\ref{T:herman} with 
$X_i:=(Y_i,\dotsc,Y_i)\in\ssl_2(\mathbb R)^n\subseteq\mathfrak X(T^1)^n\subseteq\mathfrak X(T^n)$, 
$i=1,2$, and $X_3$ so that $\exp(X_3)=R_\gamma$.
\end{proof}

\section{The exponential law}\label{S:exp}

If $\mathcal F$ is a smooth foliation of $M$ we let
$\Diff^\infty_c(M;\mathcal F)$ denote the group of compactly
supported diffeomorphisms preserving the leaves of $\mathcal F$.
This is a regular Lie group modelled on the convenient vector
space of compactly supported smooth vector fields tangential to
$\mathcal F$. The group of foliation preserving diffeomorphisms
has been studied in \cite{R95}.

\begin{lemma}\label{L:1}
Suppose $M_1$ and $M_2$ are two finite dimensional smooth manifolds and set $M:=M_1\times M_2$. Let
$\mathcal F_1$ and $\mathcal F_2$ denote the foliations with leaves $M_1\times\{\pt\}$ and 
$\{\pt\}\times M_2$ on $M$, respectively. Then the smooth map
$$
F\colon \Diff^\infty_c(M;\mathcal F_1)\times\Diff^\infty_c(M;\mathcal F_2)\to\Diff_c^\infty(M),
\quad F(g_1,g_2):=g_1\circ g_2,
$$
is a local diffeomorphism at the identity.
\end{lemma}

\begin{proof}
We proceed by constructing an inverse. To this end
let $\pi\colon M\to M_2$ denote the canonical projection and consider the smooth map
$$
\sigma_2\colon \Diff^\infty_c(M)\to C^\infty_c(M,M),\qquad\sigma_2(g)(x_1,x_2):=(x_1,\pi(g(x_1,x_2))).
$$
Since $\Diff^\infty_c(M)$ is $c^\infty$-open in $C_c^\infty(M,M)$, see \cite[Section~43.1]{KM97}, the set $\mathcal U:=\sigma_2^{-1}(\Diff^\infty_c(M))$
is a $c^\infty$-open neighborhood of the identity in $\Diff^\infty_c(M)$. Moreover, $\sigma_2$ restricts to a smooth map
$\sigma_2\colon \mathcal U\to\Diff^\infty_c(M;\mathcal F_2)$. Note that the diffeomorphism $\sigma_1(g):=g\circ\sigma_2(g)^{-1}$ preserves the leaves of the
foliation $\mathcal F_1$, and
this provides a smooth map $\sigma_1\colon \mathcal U\to\Diff^\infty_c(M;\mathcal F_1)$. We thus obtain a smooth map
$$
\sigma:=(\sigma_1,\sigma_2)\colon \mathcal U\to\Diff^\infty_c(M;\mathcal F_1)\times\Diff^\infty_c(M;\mathcal F_2)
$$
such that $F\circ\sigma=\id_{\mathcal U}$ and $\sigma(\id)=(\id,\id)$.
Since $F$ is injective, we have $\sigma(\mathcal U)=F^{-1}(\mathcal U)$, hence $\mathcal V:=\sigma(\mathcal U)$ is
$c^\infty$-open too. Clearly, $F$ restricts to a bijection, and hence diffeomorphism, $F\colon \mathcal V\cong\mathcal U$, with inverse $\sigma$.
\end{proof}

\begin{lemma}\label{L:2}
Suppose $B$ and $T$ are finite dimensional smooth manifolds, assume $T$ compact, and let $\mathcal F$ 
denote the foliation with leaves $\{\pt\}\times T$ on $B\times T$. Then the canonical bijection
$$
C^\infty_c(B,\Diff^\infty(T))\xrightarrow\cong\Diff^\infty_c(B\times T;\mathcal F)
$$
is an isomorphism of regular Lie groups.
\end{lemma}

\begin{proof}
This follows from the exponential law, see \cite[Section~42.14]{KM97}.
\end{proof}

\begin{lemma}\label{L:3}
Let $B$ be a precompact open subset in a finite dimensional smooth manifold $M$.
Then there exist compactly supported smooth vector fields $X_1,X_2,X_3$ on $M\times T^n$, tangential to the
foliation $\mathcal F$ with leaves $\{\pt\}\times T^n$, so that the map
$$
\Diff^\infty_c(B\times T^n;\mathcal F)^3\to\Diff^\infty_c(B\times T^n;\mathcal F)
$$
$$
(g_1,g_2,g_3)\mapsto[g_1,\exp(X_1)]\circ[g_2,\exp(X_2)]\circ[g_3,\exp(X_3)]
$$
admits a smooth local right inverse at the identity. Moreover, the vector fields $X_i$ may be chosen
arbitrarily close to zero with respect to the strong Whitney $C^\infty$-topology.
\end{lemma}

\begin{proof}
According to Proposition~\ref{P:herman} there exist smooth vector fields $Y_1,Y_2,Y_3$ on $T^n$ so that the map
$\Diff^\infty(T^n)^3\to\Diff^\infty(T^n)$,
$$
(h_1,h_2,h_3)\mapsto[h_1,\exp(Y_1)]\circ[h_2,\exp(Y_2)]\circ[h_3,\exp(Y_3)],
$$
admits a smooth local right inverse $\rho=(\rho_1,\rho_2,\rho_3)\colon \mathcal V\to\Diff^\infty(T^n)^3$ at the identity.
We extend the vector fields $Y_i$ in a constant manner to smooth vector fields $Z_i$ on $M\times T^n$, tangential to $\mathcal F$.
Multiplying $Z_i$ with a compactly supported smooth function which equals $1$ on $\bar B\times T^n$, we obtain
compactly supported smooth vector fields $X_1,X_2,X_3$ on $M\times T^n$, tangential to $\mathcal F$, so that
$X_i|_B=Z_i|_B$. Moreover, $X_i$ depends continuously on $Y_i$ with respect to the Whitney $C^\infty$-topologies.
Consequently, $X_i$ may be assumed arbitrarily close to zero with respect to the Whitney $C^\infty$-topology on $M$.
Consider the $c^\infty$-open neighborhood
$$
\mathcal W:=\bigl\{f\in C^\infty_c(B,\Diff^\infty(T^n))\bigm|f(B)\subseteq\mathcal V\bigr\}
$$
of the identity, and observe that the maps $(\rho_i)_*\colon \mathcal W\to C^\infty_c(B,\Diff^\infty(T^n))$
are smooth. Via the isomorphism in Lemma~\ref{L:2}, the set $\mathcal W$ corresponds to a $c^\infty$-open neighborhood $\mathcal U$
of the identity in $\Diff_c^\infty(B\times T^n;\mathcal F)$, and the maps $(\rho_i)_*$ provide a smooth map
$\sigma=(\sigma_1,\sigma_2,\sigma_3)\colon \mathcal U\to\Diff_c^\infty(B\times T^n;\mathcal F)^3$.
By construction,
$$
g=[\sigma_1(g),\exp(X_1)]\circ[\sigma_2(g),\exp(X_2)]\circ[\sigma_3(g),\exp(X_3)],\quad\text{for all $g\in\mathcal U$.}
$$
Hence $\sigma$ is the desired smooth local right inverse.
\end{proof}

\begin{lemma}\label{L:4}
Suppose $p\geq1$, $q\geq0$, set $n:=p+q$, and let $\mathcal F$ denote the foliation with leaves $\{\pt\}\times\mathbb R^q$ on $\mathbb R^n$.
Moreover, let $B$ be a precompact open subset in $\mathbb R^n$.
Then there exist compactly supported smooth vector fields $X_1,X_2,X_3$ on $\mathbb R^n$, a $c^\infty$-open neighborhood
$\mathcal U$ of the identity in $\Diff^\infty_c(B;\mathcal F)$ and smooth maps
$\sigma_1,\sigma_2,\sigma_3\colon \mathcal U\to\Diff^\infty_c(\mathbb R^n)$ so that $\sigma_i(\id)=\id$ and, for all $g\in\mathcal U$,
$$
g=[\sigma_1(g),\exp(X_1)]\circ[\sigma_2(g),\exp(X_2)]\circ[\sigma_3(g),\exp(X_3)].
$$
Moreover, the vector fields $X_i$ may be chosen arbitrarily close to zero with respect to the strong
Whitney $C^\infty$-topology on $\mathbb R^n$.
\end{lemma}

\begin{proof}
Since $p\geq1$, there exists a smooth embedding $\varphi\colon \mathbb R^p\times T^q\to\mathbb R^n$. Set
$U:=\varphi(B^p\times T^q)$, where $B^p$ denotes the open unit ball in $\mathbb R^p$. Clearly, we may assume $\bar B\subseteq U$.
Let $\mathcal G$ denote the foliation on $U$ corresponding to the foliation with leaves $\{\pt\}\times T^q$ on $B^p\times T^q$,
via the embedding $\varphi$. Furthermore, we may assume $\mathcal G|_B=\mathcal F|_B$.
According to Lemma~\ref{L:3} there exist compactly supported smooth vector fields $X_1,X_2,X_3$ on the image of $\varphi$
so that the smooth map $\Diff^\infty_c(U;\mathcal G)^3\to\Diff^\infty_c(U;\mathcal G)$,
$$
(g_1,g_2,g_3)\mapsto[g_1,\exp(X_1)]\circ[g_2,\exp(X_2)]\circ[g_3,\exp(X_3)],
$$
admits a smooth local right inverse $\rho=(\rho_1,\rho_2,\rho_3)\colon \mathcal V\to\Diff^\infty_c(U;\mathcal G)$ at the identity.
We extend the vector fields $X_i$ by zero to compactly supported smooth vector fields on $\mathbb R^n$, and observe that
these extensions may be assumed arbitrarily close to zero with respect to the Whitney $C^\infty$-topology on $\mathbb R^n$.
Let $\mathcal U$ denote the preimage of $\mathcal V$ under the canonical inclusion $\Diff^\infty_c(B;\mathcal F)\to\Diff^\infty_c(U;\mathcal G)$.
Restricting $\rho_i$, we obtain smooth maps $\sigma_1,\sigma_2,\sigma_3\colon \mathcal U\to\Diff^\infty_c(U;\mathcal G)\subseteq\Diff^\infty_c(\mathbb R^n)$
with the desired property.
\end{proof}

\begin{proposition}\label{P:exp}
Suppose $n\geq2$, and let $B$ denote a precompact open subset of $\mathbb R^n$. Then there exist compactly supported smooth vector fields $X_1,\dotsc,X_6$ on $\mathbb R^n$,
a $c^\infty$-open neighborhood $\mathcal U$ of the identity in $\Diff^\infty_c(B)$, and smooth maps
$\sigma_1,\dotsc,\sigma_6\colon \mathcal U\to\Diff^\infty_c(\mathbb R^n)$ so that $\sigma_i(\id)=\id$ and, for all $g\in\mathcal U$,
$$
g=[\sigma_1(g),\exp(X_1)]\circ\cdots\circ[\sigma_6(g),\exp(X_6)].
$$
Moreover, the vector fields $X_i$ may be chosen arbitrarily close to zero with respect to the strong
Whitney $C^\infty$-topology.
\end{proposition}

\begin{proof}
Fix $p,q\geq1$ so that $n=p+q$. Without loss of generality we may assume $B=B^p\times B^q$.
Let $\mathcal F_1$ and $\mathcal F_2$ denote the foliations with leaves $\mathbb R^p\times\{\pt\}$ and $\{\pt\}\times\mathbb R^q$
on $\mathbb R^n$, respectively. According to Lemma~\ref{L:4}, there exist compactly supported smooth vector fields
$X_4,X_5,X_6$ on $\mathbb R^n$, a $c^\infty$-open neighborhood $\mathcal W$ of the identity in $\Diff^\infty_c(B;\mathcal F_2)$,
and smooth maps $\tilde\sigma_4,\tilde\sigma_5,\tilde\sigma_6\colon \mathcal W\to\Diff^\infty_c(\mathbb R^n)$ so that $\tilde\sigma_i(\id)=\id$ and, for all $g\in\mathcal W$,
$$
g=[\tilde\sigma_4(g),\exp(X_4)]\circ[\tilde\sigma_5(g),\exp(X_5)]\circ[\tilde\sigma_6(g),\exp(X_6)].
$$
Similarly, there exist compactly supported smooth vector fields
$X_1,X_2,X_3$ on $\mathbb R^n$, a $c^\infty$-open neighborhood $\mathcal V$ of the identity in $\Diff^\infty_c(B;\mathcal F_1)$,
and smooth maps $\tilde\sigma_1,\tilde\sigma_2,\tilde\sigma_3\colon \mathcal V\to\Diff^\infty_c(\mathbb R^n)$ so that $\tilde\sigma_i(\id)=\id$ and, for all $g\in\mathcal V$,
$$
g=[\tilde\sigma_1(g),\exp(X_1)]\circ[\tilde\sigma_2(g),\exp(X_2)]\circ[\tilde\sigma_3(g),\exp(X_3)].
$$
In view of Lemma~\ref{L:1}, the composition $\Diff^\infty_c(B;\mathcal F_1)\times\Diff^\infty_c(B;\mathcal F_2)\to\Diff^\infty_c(\mathbb R^n)$
admits a smooth local right inverse $\rho=(\rho_1,\rho_2)\colon \mathcal U\to\mathcal V\times\mathcal W$ at the identity.
Hence the smooth maps $\sigma_i:=\tilde\sigma_i\circ\rho_1$, $1\leq i\leq 3$, and $\sigma_i:=\tilde\sigma_{i-3}\circ\rho_2$, $4\leq i\leq6$,
will have the desired property.
\end{proof}

\section{Smooth fragmentation}\label{S:frag}

Suppose $U\subseteq M$ is an open subset. Every compactly supported diffeomorphism of $U$ can
be regarded as a compactly supported diffeomorphism of $M$ by extending it identically outside $U$.
The resulting injective homomorphism $\Diff_c^\infty(U)\to\Diff^\infty_c(M)$ is clearly smooth. Note, however, that
a curve in $\Diff^\infty_c(U)$, which is smooth when considered as a curve in $\Diff_c^\infty(M)$,
need not be smooth as a curve into $\Diff^\infty_c(U)$. Nevertheless, if there exists a closed subset $A$ of $M$ with
$A\subseteq U$ and if the curve has support contained in $A$, then one can conclude that the curve is smooth in
$\Diff^\infty_c(U)$, too. This follows immediately from the characterization of smooth curves given at the beginning of Section~\ref{S:intro}.

The following is a folklore statement which, in one form or the other, can be found all over the literature on
diffeomorphism groups \cite{B97}. For the reader's convenience we include a version emphasizing the
smoothness of the construction.

\begin{proposition}[Fragmentation]\label{P:frag}
Let $M$ be a smooth manifold of dimension $n$, and
suppose $U_1,\dotsc,U_k$ is an open covering of $M$, ie.\ $M=U_1\cup\cdots\cup U_k$.
Then the smooth map
$$
P\colon \Diff^\infty_c(U_1)\times\cdots\times\Diff^\infty_c(U_k)
\to\Diff^\infty_c(M),
\quad
P(g_1,\dotsc,g_k):=g_1\circ\cdots\circ g_k,
$$
admits a smooth local right inverse at the identity.
\end{proposition}

\begin{proof}
Let $\pi\colon TM\to M$ denote the projection of the tangent bundle.
Choose a Riemannian metric on $M$ and let $\exp$ denote the
corresponding exponential map. Choose an open neighborhood of the
zero section $W\subseteq TM$ such that
$(\pi,\exp)\colon W\to M\times M$ is a diffeomorphism onto an open
neighborhood of the diagonal. Fix a linear connection $\nabla$ on $TM$, let
$W'\subseteq T^*M\otimes TM$ be an open neighborhood
of the zero section, and set
$$
\mathcal W:=\bigl\{X\in\mathfrak X_c(M)\bigm|\img(X)\subseteq W,\, \img(\nabla X)\subseteq W'\bigr\}.
$$
This is a $C^1$-open, hence $c^\infty$-open, zero neighborhood in $\mathfrak X_c(M)$.
For $X\in\mathcal W$ define $f_X:=\exp\circ X\in C^\infty_c(M,M)$.
Choosing $W$ and $W'$ sufficiently small, we may assume that every $f_X$, $X\in\mathcal W$, is a
diffeomorphism \cite{H76}. The map
\begin{equation}\label{E:chart}
\mathcal W\to\Diff^\infty_c(M),\quad X\mapsto f_X
\end{equation}
provides a chart of $\Diff^\infty_c(M)$ centered
at the identity. This is the standard way to put a smooth structure on
$\Diff^\infty_c(M)$, see \cite[Section~42.1]{KM97}.

Choose a smooth partition of unity $\lambda_1,\dotsc,\lambda_k$ with
$\supp(\lambda_i)\subseteq U_i$, $1\leq i\leq k$, and define
$$
\mathcal V:=\bigl\{X\in\mathcal W\bigm|\text{$(\lambda_1+\cdots+\lambda_i)X\in\mathcal W$ for all $1\leq i\leq k$}\bigr\}.
$$
This is a $C^1$-open, hence $c^\infty$-open, zero neighborhood in $\mathfrak X_c(M)$, and $\mathcal V\subseteq\mathcal W$.
For $X\in\mathcal V$ and $1\leq i\leq k$ set $X_i:=(\lambda_1+\cdots+\lambda_i)X$, and note that
$f_{X_i}\in\Diff^\infty_c(M)$. Clearly, $\supp(f_{X_1})\subseteq\supp(\lambda_1)\subseteq U_1$.
Moreover, for $1<i\leq k$ we have $X_{i-1}=X_i$ on $M\setminus\supp(\lambda_i)$, and thus
$f_{X_{i-1}}=f_{X_i}$ on $M\setminus\supp(\lambda_i)$.
We conclude that the diffeomorphism $(f_{X_{i-1}})^{-1}\circ f_{X_i}$ has compact support contained in
$\supp(\lambda_i)\subseteq U_i$, $1<i\leq k$.

Let $\mathcal U\subseteq\Diff_c^\infty(M)$ denote the $c^\infty$-open neighborhood
of the identity in $\Diff^\infty_c(M)$ corresponding to $\mathcal V$ via \eqref{E:chart}, and define a map
$$
\sigma\colon \mathcal U\to\Diff^\infty_c(U_1)\times\cdots\times\Diff^\infty_c(U_k)
$$
$$
\sigma(f_X):=\Bigl(f_{X_1},(f_{X_1})^{-1}\circ f_{X_2},
(f_{X_2})^{-1}\circ f_{X_3},
\dotsc,(f_{X_{k-1}})^{-1}\circ f_{X_k}\Bigr).
$$
Since the support of $(f_{X_{i-1}})^{-1}\circ f_{X_i}$ is contained in $\supp(\lambda_i)\subseteq U_i$, it is clear
that $\sigma$ is a smooth map, cf.\ the remark at the beginning of this section. Obviously we have
$P(\sigma(f_X))=f_{X_k}=f_X$ and thus $P\circ\sigma=\id_{\mathcal U}$.
Moreover, it is immediate from the construction that $\sigma(\id)=\sigma(f_0)=(\id,\dotsc,\id)$.
\end{proof}

Combining Propositions~\ref{P:exp} and \ref{P:frag} permits to show the following result.

\begin{proposition}\label{P:main}
Let $M$ be a smooth manifold of dimension $n\geq2$, set $N:=6(n+1)$, and
suppose $U$ and $V$ are open subsets of $M$ such that $\bar U\subseteq V$.
Then there exist smooth complete vector fields
$X_1,\dotsc,X_N$ on $M$ with $\supp(X_i)\subseteq V$, a $c^\infty$-open neighborhood $\mathcal U$ of the identity in
$\Diff^\infty_c(U)$ and smooth maps $\sigma_1,\dotsc,\sigma_N\colon\mathcal U\to\Diff^\infty_c(V)$ such that
$\sigma_i(\id)=\id$ and, for all $g\in\mathcal U$,
$$
g=[\sigma_1(g),\exp(X_1)]\circ\cdots\circ[\sigma_N(g),\exp(X_N)].
$$
Moreover, the vector fields $X_i$ may be chosen arbitrarily close to zero with respect to the strong
Whitney $C^\infty$-topology.
\end{proposition}

\begin{proof}
It is well known that there exist open subsets $U_1,\dotsc,U_{n+1}$ of $M$ so that
$$
\bar U\subseteq U_1\cup\cdots\cup U_{n+1}\subseteq V
$$
and such that each $U_i$, $1\leq i\leq n+1$, is diffeomorphic to a disjoint union of copies of the open unit ball $B^n$.
Moreover, we may assume that there exist embeddings
$\varphi_i\colon \bigsqcup_{\alpha\in A_i}\mathbb R^n\to V$, with index sets $A_i$, so that
$$
\textstyle
\varphi_i\bigl(\,\bigsqcup_{\alpha\in A_i}B^n\bigr)=U_i,\qquad1\leq i\leq n+1.
$$
Applying Proposition~\ref{P:exp} to each connected component of $U_i$, we find complete vector fields
$X_{i,1},\dotsc,X_{i,6}$ on $M$ with $\supp(X_{i,j})\subseteq V$, a $c^\infty$-open neighborhood $\mathcal U_i$ of the identity in $\Diff_c^\infty(U_i)$
and smooth maps $\sigma_{i,1},\dotsc,\sigma_{i,6}\colon \mathcal U_i\to\Diff^\infty_c(U_i)\subseteq\Diff_c^\infty(V)$ so that
$\sigma_{i,j}(\id)=\id$ and, for all $g\in\mathcal U_i$,
$$
g=[\sigma_{i,1}(g),\exp(X_{i,1})]\circ\cdots\circ[\sigma_{i,6}(g),\exp(X_{i,6})].
$$
Moreover, the vector fields $X_{i,j}$ may be chosen arbitrarily close to zero with respect to the strong Whitney $C^\infty$-topology on $M$.
In view of Proposition~\ref{P:frag}, the map
\begin{align*}
\Diff_c^\infty(U\cap U_1)\times\cdots\times\Diff_c^\infty(U\cap U_{n+1})&\to\Diff_c^\infty(U),
\\
(g_1,\dotsc,g_{n+1})&\mapsto g_1\circ\cdots\circ g_{n+1}
\end{align*}
admits a local smooth right inverse at the identity. Combining this with the $\sigma_{i,j}$ above,
we immediately obtain the statement.
\end{proof}

Specializing Proposition~\ref{P:main} to $U=M$, we obtain Theorem~\ref{T:mainp}.

\section{Reducing the number of commutators}\label{S:proof}

Proceeding as in \cite{BIP08} permits to reduce the number of commutators considerably, see
also \cite{T08} and \cite{T09}.

\begin{proposition}\label{P:rho}
Let $M$ be a smooth manifold of dimension $n\geq2$ and put $N=6(n+1)$.
Moreover, let $U$ an open subset of $M$ and suppose $\phi\in\Diff^\infty(M)$, not necessarily with
compact support, such that the closures of the subsets
$$
U,\,\phi(U),\,\phi^2(U),\,\dotsc,\,\phi^N(U)
$$
are mutually disjoint.
Then there exists a smooth complete vector field $X$ on $M$, a $c^\infty$-open neighborhood
$\mathcal U$ of the identity in $\Diff^\infty_c(U)$, and smooth maps
$\varrho_1,\varrho_2\colon\mathcal U\to\Diff^\infty_c(M)$ so that
$\varrho_1(\id)=\varrho_2(\id)=\id$ and, for all $g\in\mathcal U$,
$$
g=[\varrho_1(g),\phi]\circ[\varrho_2(g),\exp(X)].
$$
Moreover, the vector field $X$ may be chosen arbitrarily close to zero in the strong Whitney $C^\infty$-topology on $M$.
\end{proposition}

\begin{proof}
Clearly, there exists an open subset $V$ of $M$ with $\bar U\subseteq V$ such that the open subsets
$V,\phi(V),\phi^2(V),\dotsc,\phi^N(V)$ are mutually disjoint.
Define a smooth map $\psi\colon\Diff^\infty_c(V)^N\to\Diff_c^\infty(M)$,
$$
\psi(g_1,\dotsc,g_N):=\prod_{i=1}^N\phi^{i-1}\circ g_i\circ\cdots\circ g_N\circ\phi^{-(i-1)}.
$$
If $g,h\in\Diff_c^\infty(V)$ and $0\leq i\neq j\leq N$, then the diffeomorphisms
$\phi^i\circ g\circ\phi^{-i}$ and $\phi^j\circ h\circ\phi^{-j}$ commute as their supports
are disjoint. Hence
\begin{equation}\label{E:abc}
g_1\circ\cdots\circ g_N=[\psi(g_1,\dotsc,g_N),\phi]\circ\prod_{i=1}^N\phi^i\circ g_i\circ\phi^{-i},
\end{equation}
for all $g_1,\dotsc,g_N\in\Diff_c^\infty(V)$.
According to Proposition~\ref{P:main} there
exist smooth complete vector fields $X_1,\dotsc,X_N$ on $M$ with $\supp(X_i)\subseteq V$, a $c^\infty$-open neighborhood $\mathcal U$ of the identity in $\Diff^\infty_c(U)$
and smooth maps $\sigma_1,\dotsc,\sigma_N\colon\mathcal U\to\Diff_c^\infty(V)$ so that $\sigma_i(\id)=\id$ and,
for all $g\in\mathcal U$,
$$
g=[\sigma_1(g),\exp(X_1)]\circ\cdots\circ[\sigma_N(g),\exp(X_N)].
$$
Combining this with \eqref{E:abc} we obtain, for all $g\in\mathcal U$,
$$
g=[\varrho_1(g),\phi]\circ[\varrho_2(g),\exp(X)],
$$
where $\varrho_1,\varrho_2\colon\mathcal U\to\Diff_c^\infty(M)$,
\begin{align*}
\varrho_1(g)&:=\psi\bigl([\sigma_1(g),\exp(X_1)],\dotsc,[\sigma_N(g),\exp(X_N)]\bigr),
\\
\varrho_2(g)&:=\prod_{i=1}^N\phi^i\circ\sigma_i(g)\circ\phi^{-i},
\end{align*}
and $X:=\sum_{i=1}^N\phi^i_*X_i$, i.e.\ $\exp(X)=\prod_{i=1}^N\phi^i\circ\exp(X_i)\circ\phi^{-i}$.
Since $X$ depends continuously on $X_1,\dotsc,X_N$,
we may assume $X$ to be arbitrarily close to zero with respect to the strong Whitney $C^\infty$-topology.
Clearly, $\varrho_1$ and $\varrho_2$ are smooth, and we have $\varrho_1(\id)=\varrho_2(\id)=\id$.
\end{proof}

\begin{lemma}\label{L:displ}
Let $M$ be a smooth manifold of dimension $n$.
Then there exists an open covering $M=U_1\cup U_2\cup U_3$ and
smooth complete vector fields $X_1,X_2,X_3$ on $M$ so that
$\exp(X_1)(U_1)\subseteq U_2$, $\exp(X_2)(U_2)\subseteq U_3$, and such that
the closures of the sets
$$
U_3,\,\,\,\exp(X_3)(U_3),\,\,\,\exp(X_3)^2(U_3),\,\,\dotsc
$$
are mutually disjoint. Moreover, the vector fields $X_1,X_2,X_3$ may be chosen arbitrarily close to zero
with respect to the strong Whitney $C^0$-topology.
If $M$ admits a proper (circle valued) Morse function whose critical points all have index $0$ or $n$,
then we may, moreover, choose $U_1=\emptyset$ and $X_1=0$.
\end{lemma}

\begin{proof}
Fix a proper Morse function $f\colon M\to\R$, see \cite{H76}, and let
$\mathcal X$ denote the set of critical points of $f$. Since $\mathcal X$ is a discrete subset of $M$,
the properness of $f$ implies that the critical values form a discrete subset of $\R$.
We may, moreover, assume that each critical level contains precisely one critical point, see \cite{Milnor}.
For notational brevity we will write $\ind(t)$ for the Morse index of the unique critical point
corresponding to a critical value $t$ of $f$. For any subset $X\subseteq\R$, we introduce the notation
$M_X:=f^{-1}(X)\subseteq M$.

Let $\mathcal U$ be a zero neighborhood in the space of vector fields on $M$ with respect to the strong
Whitney $C^0$-topology. A neighborhood basis in this topology is obtained by taking neighborhoods of the
zero section in $TM$ and considering vector fields which take values in this subset.
We will assume that $\mathcal U$ is of this form, for a neighborhood of
the zero section which is fiberwise convex. Particularly, $\mathcal U$ is invariant with respect to
multiplication by functions whose modulus does not exceed $1$. Below we will show that the following
assertions hold true:
\begin{enumerate}[(a)]
\item\label{I:a}
Suppose $t$ is a regular value of $f$. Then, for sufficiently small $\varepsilon>0$ and all $\eta>0$
there exists a vector field $Y_2^{\varepsilon,\eta}\in\mathcal U$ so that
$\supp(Y_2^{\varepsilon,\eta})\subseteq M_{(t,t+\varepsilon+\eta)}$ and
$$
\exp(Y_2^{\varepsilon,\eta})(M_{(t,t+\varepsilon)})\subseteq M_{(t,t+\eta)}.
$$
\item\label{I:b}
Suppose $t$ is a regular value of $f$. Then, for sufficiently small $\eta>0$
there exists a vector field $Y_3^\eta\in\mathcal U$ with $\supp(Y_3^\eta)\subseteq M_{(t-2\eta,t+2\eta)}$ such that the closures of the subsets
$$
M_{(t-\eta,t+\eta)},\,\,\exp(Y_3^\eta)\bigl(M_{(t-\eta,t+\eta)}\bigr),\,\,\exp(Y_3^\eta)^2\bigl(M_{(t-\eta,t+\eta)}\bigr),\,\,\dotsc
$$
are mutually disjoint.
\item\label{I:c}
Suppose $t$ is a critical value of $f$ and $0<\ind(t)<n$.
Then, for sufficiently small $\varepsilon>0$ and all $\eta>0$, there exists an open covering
$M_{(t-\varepsilon,t+\varepsilon)}=V_1^\varepsilon\cup V_2^\varepsilon$ and vector fields
$Y_1^{\varepsilon,\eta},Y_2^{\varepsilon,\eta}\in\mathcal U$ such that
$$
\exp(Y_1^{\varepsilon,\eta})(V_1^\varepsilon)\subseteq V_2^\varepsilon,
$$
$\supp(Y_1^{\varepsilon,\eta})\subseteq M_{(t-\varepsilon,t+\varepsilon+\eta)}$,
$\supp(Y_2^{\varepsilon,\eta})\subseteq M_{(t-\varepsilon,t+\varepsilon+\eta)}$, and
$$
\exp(Y_2^{\varepsilon,\eta})(V_2^\varepsilon)\subseteq M_{(t-\varepsilon,t-\varepsilon+\eta)}.
$$
\item\label{I:d}
Suppose $t$ is a critical value of $f$ with $\ind(t)=n$, i.e.\ a local maximum.
Then, for sufficiently small $\varepsilon>0$ and all $\eta>0$, there exists a vector field
$Y_2^{\varepsilon,\eta}\in\mathcal U$ such that $\supp(Y_2^{\varepsilon,\eta})\subseteq M_{(t-\varepsilon,t+\varepsilon+\eta)}$ and
$$
\exp(Y_2^{\varepsilon,\eta})(M_{(t-\varepsilon,t+\varepsilon)})\subseteq M_{(t-\varepsilon,t-\varepsilon+\eta)}.
$$
\item\label{I:e}
Suppose $t$ is a critical value of $f$ with $\ind(t)=0$, i.e.\ a local minimum.
Then, for sufficiently small $\varepsilon>0$ and all $\eta>0$, there exists a vector field
$Y_2^{\varepsilon,\eta}\in\mathcal U$ such that $\supp(Y_2^{\varepsilon,\eta})\subseteq M_{(t-\varepsilon-\eta,t+\varepsilon)}$ and
$$
\exp(Y_2^{\varepsilon,\eta})(M_{(t-\varepsilon,t+\varepsilon)})\subseteq M_{(t+\varepsilon-\eta,t+\varepsilon)}.
$$
\end{enumerate}

We postpone the proof of these statements and first indicate how the lemma can be derived from them.
For each critical value $t$ of $f$ we choose $\varepsilon_t>0$ and $\eta_t>0$ as in (\ref{I:c}), (\ref{I:d}) or (\ref{I:e}), respectively, depending on $\ind(t)$.
Shrinking $\varepsilon_t$ and $\eta_t$ we may assume that the intervals $[t-\varepsilon_t-\eta_t,t+\varepsilon_t+\eta_t]$ are mutually disjoint.
We let
\begin{align*}
U_1&:=\bigsqcup_{0<\ind(t)<n}V_1^{\varepsilon_t}
&X_1&:=\sum_{0<\ind(t)<n}Y^{\varepsilon_t,\eta_t}_1
\\
U_2'&:=\bigsqcup_{0<\ind(t)<n}V_2^{\varepsilon_t}
&X_2'&:=\sum_{0<\ind(t)<n}Y^{\varepsilon_t,\eta_t}_2
\end{align*}
where $t$ runs through all critical values of $f$ corresponding to critical points of index different from $0$ and $n$.
Moreover, $V_1^{\varepsilon_t}$, $V_2^{\varepsilon_t}$, $Y^{\varepsilon_t,\eta_t}_1$, $Y^{\varepsilon_t,\eta_t}_2$,
denote the open subsets and vector fields constructed in (\ref{I:c}), respectively.
Furthermore, we set
\begin{align*}
U_2''&:=U_2'\sqcup\bigsqcup_{\ind(t)=0,n}M_{(t-\varepsilon_t,t+\varepsilon_t)}
&X_2''&:=X_2'+\sum_{\ind(t)=0,n}Y_2^{\varepsilon_t,\eta_t},
\end{align*}
where $t$ runs through all critical values $t$ of $f$ corresponding to critical points of index $0$ or $n$,
and the vector fields $Y_2^{\varepsilon_t,\eta_t}$ are the ones constructed in (\ref{I:d}) or (\ref{I:e}), respectively.
By construction we have $X_1,X_2''\in\mathcal U$, $\exp(X_1)(U_1)\subseteq U_2''$, $U_1\cup U_2''=\bigsqcup_tM_{(t-\varepsilon_t,t+\varepsilon_t)}$, where the latter
union is over all critical values, and
$$
\exp(X_2'')(U_2'')\subseteq\bigsqcup_{\ind(t)>0}M_{(t-\varepsilon_t,t-\varepsilon_t+\eta_t)}\sqcup\bigsqcup_{\ind(t)=0}M_{(t+\varepsilon_t-\eta_t,t+\varepsilon_t)}.
$$
Note that we are still free to shrink all $\eta_t$ without affecting these properties.
Using (\ref{I:a}) and (\ref{I:b}) it is now clear how to complete the construction,
the open sets $U_2$ and $U_3$ can be chosen to be of the form $U_2=U_2''\sqcup M_I$ and $U_3=M_J$ where $I$ and $J$
are disjoint unions of suitably chosen intervals.

It thus remains to verify assertions (\ref{I:a}) through (\ref{I:e}) above.
If $t$ is a regular value of $f$, then there exists an open interval $I$ containing $t$ and a diffeomorphisms
$M_I\cong M_t\times I$ intertwining the map $f\colon M_I\to I$ with the standard projection $M_t\times I\to I$, see \cite{Milnor}.
In order to prove statements (\ref{I:a}) and (\ref{I:b}) it thus suffices to write down appropriate vector
fields on $I$ which is straightforward.

Now suppose $t$ is a critical value of $f$ with corresponding critical point $y$ of index $k$.
Choose a Morse chart \cite{Milnor} centered at $y$, i.e.\ $(x_1,\dotsc,x_n)\colon W\to\R^n$
are local coordinates such that
$$
f|_W=f(y)-x_1^2-\cdots-x_k^2+x_{k+1}^2+\cdots+x_n^2
$$
on the open neighborhood $W$ of $y$. We fix a Riemannian metric $g$ on $M$ such that its restriction to $W$
is the standard Euclidean metric, i.e.\
$$
g|_W=dx^1\otimes dx^1+\cdots+dx^n\otimes dx^n.
$$
The gradient vector field $X:=-\grad(f)$ thus has a simple linear form on $W$,
$$
X|_W=2x_1\partial_1+\cdots+2x_k\partial_k-2x_{k+1}\partial_{k+1}-\cdots-2x_n\partial_n,
$$
where $\partial_i=\frac\partial{\partial x_i}$, $1\leq i\leq n$, denote the associated coordinate vector fields. Set
$$
V_1^\varepsilon:=M_{(t-\varepsilon,t+\varepsilon)}\cap\{x\in W:x_1^2+\cdots+x_k^2<\varepsilon/2\}
$$
and
$$
V_2^\varepsilon:=M_{(t-\varepsilon,t+\varepsilon)}\setminus\{x\in W:x_1^2+\cdots+x_k^2\leq\varepsilon/3\}.
$$
Clearly, $M_{(t-\varepsilon,t+\varepsilon)}=V_1^\varepsilon\cup V_2^\varepsilon$, and $V_1^\varepsilon$ is open.
For sufficiently small $\varepsilon>0$ the set $V_2^\varepsilon$ will be open too.
By construction, the flow lines of $X$ starting at points in $V_2^\varepsilon$ remain bounded away from the critical point $y$ by a positive distance.
For sufficiently small $\varepsilon>0$, we thus have $\exp(X)(V_2^\varepsilon)\subseteq M_{(-\infty,t-\varepsilon)}$.
Multiplying $X$ with an appropriate function, we get a vector field $Y_2^{\varepsilon,\eta}\in\mathcal U$ such that
$\supp(Y_2^{\varepsilon,\eta})\subseteq M_{(t-\varepsilon,t+\varepsilon+\eta)}$, and
$$
\exp(Y_2^{\varepsilon,\eta})(V^\varepsilon_2)\subseteq M_{(t-\varepsilon,t-\varepsilon+\eta)}.
$$
Moreover, if $k\neq0$, then after possibly shrinking $\varepsilon$, it is straightforward to write down a vector field $Y^{\varepsilon,\eta}_1\in\mathcal U$, supported in
$W$, such that $\exp(Y_1^{\varepsilon,\eta})(V_1^\varepsilon)\subseteq V_2^\varepsilon$ and $\supp(Y_1^\varepsilon)\subseteq M_{(t-\varepsilon,t+\varepsilon+\eta)}$.
This shows (\ref{I:c}).

To see (\ref{I:d}), i.e.\ the case $k=n$, note first that $M_{(t-2\varepsilon,t+2\varepsilon)}\cap W$ is a connected component of $M_{(t-2\varepsilon,t+2\varepsilon)}$,
for sufficiently small $\varepsilon>0$. On $M_{(t-\varepsilon,t+\varepsilon+\eta)}\setminus W$ we construct the vector field $Y_2^{\varepsilon,\eta}$
as above. On $M_{(t-2\varepsilon,t+2\varepsilon)}\cap W$ it is straightforward to write down a vector field $Y_2^{\varepsilon,\eta}$ with the desired properties explicitly.

Finally, statement (\ref{I:e}) follows from (\ref{I:d}) by considering $-f$.
\end{proof}

We are now in a position to complete the proof of Theorem~\ref{T:main}.
Fix an open covering $M=U_1\cup U_2\cup U_3$ and smooth complete vector fields $X_1,X_2,X_3$
as in Lemma~\ref{L:displ} above.
According to Proposition~\ref{P:frag} there exists a $c^\infty$-open neighborhood $\mathcal V$
of the identity in $\Diff^\infty_c(M)$ and smooth maps $\tilde\sigma_i\colon\mathcal V\to\Diff_c^\infty(U_i)$
such that $\tilde\sigma_i(\id)=\id$ and, for all $g\in\mathcal V$,
$$
g=\tilde\sigma_1(g)\circ\tilde\sigma_2(g)\circ\tilde\sigma_3(g).
$$
A trivial computation shows that for all $g\in\mathcal V$ we have
$$
g=[\sigma_1(g),\exp(X_1)]\circ[\sigma_2(g),\exp(X_2)]\circ\phi(g),
$$
where $\sigma_1(g):=\tilde\sigma_1(g)$ and
\begin{align*}
\sigma_2(g)&:=\exp(X_1)\circ\sigma_1(g)\circ\exp(X_1)^{-1}\circ\tilde\sigma_2(g),
\\
\phi(g)&:=\exp(X_2)\circ\sigma_2(g)\circ\exp(X_2)^{-1}\circ\tilde\sigma_3(g).
\end{align*}
Clearly, these expressions define smooth maps $\sigma_1\colon\mathcal V\to\Diff_c^\infty(U_1)\subseteq\Diff_c^\infty(M)$,
$\sigma_2\colon\mathcal V\to\Diff_c^\infty(U_2)\subseteq\Diff_c^\infty(M)$, and $\phi\colon\mathcal V\to\Diff_c^\infty(U_3)$
such that $\sigma_1(\id)=\sigma_2(\id)=\phi(\id)=\id$.
In view of Proposition~\ref{P:rho} there exists a smooth complete vector field $X_4$ on $M$, a $c^\infty$-open neighborhood $\mathcal W$
of the identity in $\Diff_c^\infty(U_3)$, and smooth maps $\varrho_1,\varrho_2\colon\mathcal W\to\Diff^\infty_c(M)$ so that $\varrho_1(\id)=\varrho_2(\id)=\id$,
and
$$
h=[\varrho_1(h),\exp(X_3)]\circ[\varrho_2(h),\exp(X_4)],
$$
for all $h\in\mathcal W$.
Note that $\mathcal U:=\phi^{-1}(\mathcal W)$ is a $c^\infty$-open neighborhood of the identity in $\Diff^\infty_c(M)$.
Defining smooth maps $\sigma_3,\sigma_4\colon\mathcal U\to\Diff_c^\infty(M)$ by
$\sigma_3(g):=\varrho_1(\phi(g))$ and $\sigma_4(g)=\varrho_2(\phi(g))$, we obtain $\sigma_3(\id)=\sigma_4(\id)=\id$ and
$$
g=[\sigma_1(g),\exp(X_1)]\circ[\sigma_2(g),\exp(X_2)]\circ[\sigma_3(g),\exp(X_3)]\circ[\sigma_4(g),\exp(X_4)],
$$
for all $g\in\mathcal U$.
It follows immediately from the corresponding assertions in Proposition~\ref{P:rho} and Lemma~\ref{L:displ} that
the vector fields $X_1,\dotsc,X_4$ may be chosen arbitrarily close to zero in the strong Whitney $C^0$-topology.
If $M$ admits a (circle) valued proper Morse function whose critical points all have index $0$ or $n$,
then we may assume $X_1=0$, see Lemma~\ref{L:displ}. Whence, in this case,
$$
g=[\sigma_2(g),\exp(X_2)]\circ[\sigma_3(g),\exp(X_3)]\circ[\sigma_4(g),\exp(X_4)],
$$
for all $g\in\mathcal U$. This completes the proof of Theorem~\ref{T:main}.


\begin{thebibliography}{XX}

\bibitem{B97}
  A. Banyaga,
  \emph{The structure of classical diffeomorphism groups,}
  Kluwer Academic Publishers Group, Dordrecht, 1997.

\bibitem{BIP08}
  D. Burago, S. Ivanov and L. Polterovich,
  \emph{Conjugation-invariant norms on groups of geometric origin,}
  Groups of diffeomorphisms, 221--250,
  Adv. Stud. Pure Math. \textbf{52}, Math. Soc. Japan, Tokyo, 2008.

\bibitem{E70}
  D.B.A. Epstein,
  \emph{The simplicity of certain groups of homeomorphisms,}
  Compositio Math. \textbf{22}(1970), 165--173.

\bibitem{E84}
  D.B.A. Epstein,
  \emph{Commutators of $C^\infty$-diffeomorphisms. Appendix to: ``A curious
    remark concerning the geometric transfer map'' by John N. Mather,}
  Comment. Math. Helv. \textbf{59}(1984), 111--122.

\bibitem{HT03}
  S. Haller and J. Teichmann,
  \emph{Smooth perfectness through decomposition of
        diffeomorphisms into fiber preserving ones,}
  Ann. Global Anal. Geom. \textbf{23}(2003), 53--63.

\bibitem{H71}
  M.R. Herman,
  \emph{Simplicit\'e du groupe des diff\'eomorphismes de
        classe $C^\infty$, isotopes \`a l'identit\'e, du tore
        de dimension $n$,}
  C. R. Acad. Sci. Paris Sr. A \textbf{273}(1971), 232--234.

\bibitem{H73}
  M.R. Herman,
  \emph{Sur le groupe des difféomorphismes du tore,}
  Ann. Inst. Fourier (Grenoble) \textbf{23}(1973), 75--86.

\bibitem{H76}
  M.W. Hirsch,
  \emph{Differential topology.}
  Graduate Texts in Mathematics \textbf{33}, Springer, 1976.

\bibitem{KM97}
  A. Kriegl and P.W. Michor,
  \emph{The convenient setting of global analysis.}
  Mathematical Surveys and Monographs \textbf{53},
  American Mathematical Society, 1997.

\bibitem{M74}
  J.N. Mather,
  \emph{Commutators of diffeomorphisms,}
  Comment. Math. Helv. \textbf{49}(1974), 512--528.

\bibitem{M75}
  J.N. Mather,
  \emph{Commutators of diffeomorphisms. II,}
  Comment. Math. Helv. \textbf{50}(1975), 33--40.

\bibitem{M84}
  J.N. Mather,
  \emph{A curious remark concerning the geometric transfer map,}
  Comment. Math. Helv. \textbf{59}(1984), 86--110.

\bibitem{Milnor}
  J. Milnor,
  \emph{Morse theory.}
  Annals of Mathematics Studies \textbf{51},
  Princeton University Press, Princeton, N.J. 1963.

\bibitem{Novikov}
  S.P. Novikov,
  \emph{Multivalued functions and functionals. An analogue of the Morse theory,}
  Dokl. Akad. Nauk SSSR \textbf{260}(1981), 31--35.
  English translation: Soviet Math. Dokl. \textbf{24}(1981), 222--226(1982).

\bibitem{Pajitnov}
  A.V. Pajitnov,
  \emph{Circle-valued Morse theory.}
  de Gruyter Studies in Mathematics \textbf{32},
  Walter de Gruyter \& Co., Berlin, 2006.

\bibitem{R95}
  T. Rybicki,
  \emph{The identity component of the leaf preserving diffeomorphism group is perfect,}
  Monatsh. Math. \textbf{120}(1995), 289--305.

\bibitem{R11}
  T. Rybicki,
  \emph{Locally continuously perfect groups of homeomorphisms,}
  Ann. Glob. Anal. Geom. \textbf{40}(2011), 191--202.

\bibitem{T74}
  W. Thurston,
  \emph{Foliations and groups of diffeomorphisms,}
  Bull. Amer. Math. Soc. \textbf{80}(1974), 304--307.

\bibitem{T81}
  T. Tsuboi,
  \emph{On $2$-cycles of $B\Diff(S^1)$ which are represented by foliated $S^1$-bundles over $T^2$,}
  Ann. Inst. Fourier \textbf{31}(1981), 1--59.

\bibitem{T08}
  T. Tsuboi, 
  \emph{On the uniform perfectness of diffeomorphism groups.} 
  Groups of diffeomorphisms, 505--524, Adv. Stud. Pure Math. \textbf{52}, 
  Math. Soc. Japan, Tokyo, 2008.

\bibitem{T09}
  T. Tsuboi, 
  \emph{On the uniform simplicity of diffeomorphism groups.} 
  Differential geometry, 43--55, World Sci. Publ., Hackensack, NJ, 2009.

\bibitem{yoccoz95}
  J.-C. Yoccoz,
  \emph{Petits diviseurs en dimension $1$}, Ast\'erisque \textbf{231}(1995).
\end{thebibliography}
\end{document}